\documentclass{article}
\usepackage{amssymb}
\usepackage{amsmath}
\usepackage{graphics}
\usepackage{hyperref}
\usepackage{xcolor}
\usepackage{tikz}
\usepackage{pgfplots}
\usepackage{amsfonts}

\usetikzlibrary{arrows.meta}
\usetikzlibrary{decorations.text}
\usetikzlibrary {calc,quotes}
\pgfplotsset{compat=1.18}
\newtheorem{theorem}{Theorem}

\newtheorem{lemma}[theorem]{Lemma}

\newenvironment{proof}[1][Proof]{\noindent\textbf{#1.} }{\ \rule{0.5em}{0.5em}}

\begin{document}

\title{On a simple family of Ricci-flows of O(3)-hypersurfaces:\\
Counterexamples to the theorem on the preservation of positive definite
Ricci curvature everywhere}
\author{Wu-Chung Hsiang\thanks{%
Department of Mathematics, Princeton University, Princeton, NJ 08544 } \ and
Wu-Yi Hsiang\thanks{%
Department of Mathematics, University of California, Berkeley, CA 94720}}
\maketitle

\underline{Abstract}

This note provides a concise derivation of the equation of Ricci flows on
O(n)-rotational hypersurfaces in $\mathbb{E}^{n+1}$ that consists of a pair
of parabolic equations (cf. Lemma 2 of \S 4) via simple applications of
classical theorems of rotational surfaces in $\mathbb{E}^{3}$, governing the
Ricci-evolution of their generating plane curves. Furthermore, the second
equation can be transformed to that of the evolution equation of their
"curvature profiles" (cf. Lemmas 3 and 3' of \S 4). In \S 5, we consider a
specific family of initial generating curves of "curvature profiles" that
are only tiny modifications of the constant function $K=\frac{1}{R},R\gg 1$,
in a very small vicinity of the symmetric middle point as indicated in Fig.
2. Note that the short-time unique existence of analytic solutions of such a
simple system of parabolic equations follows simply via Nash-Moser Theory.
The simple estimate of \S 5 shows that the curvature $k_{1}(\gamma _{t}(0))$
is already becoming negative in a very short time (cf. \S 5.3).

\section{Introduction}

Historically, the formulation as well as the study of differentiable spaces
with metric structure (e.g. definition of $ds^{2}$) was first initiated in
1842 by Jacobi in his geometrization of Lagrange's least action principle of
classical mechanics, and then in 1854 by Riemann as a kind of far-reaching
abstract generalization of classical surface theory (most likely independent
of, and unaware of Jacobi's work), nowadays commonly referred to as
Riemannian manifolds and Riemannian geometry\footnote{%
Jacobi's focus was on the 3-body problem in celestial mechanics, but it was
sidetracked by the development of Hamilton-Jacobi theory and was cut short
by his early death in 1851.}. For examples, let us just mention here
Jacobi's foundational work on the global geometry of geodesics such as
Jacobi's equation and the second variational formula, and Riemann's
foundational theory on local geometric invariants on the Riemannian
curvature tensor and covariant differentiation, which is a far-reaching
generalization of Gauss' classical theory on surfaces of $\mathbb{E}^{3}$;
while the global manifestations of various simple curvature assumptions
naturally become interesting \ and/or important testing problems of global
geometry.

Inspired by the Cartan-Hadamard Theorem on simply connected manifolds $M^{n}$
with \textit{non-positive sectional curvature }everywhere and in every
direction (namely, diffeomorphic to $\mathbb{E}^{n}$ as a consequence of
Jacobi's equation), the topic of the determination of homeomorphic (resp.
diffeomorphic) type of simply connected complete $M^{n}$ with various kinds
of everywhere curvature positive conditions has been a long-standing,
prominent topic of global differential geometry for seven decades up to the
2010s.

In retrospect, the works on this topic roughly consists of two periods with
rather different approaches and geometric-analytical techniques. Briefly
speaking, the major results in this direction before the 1980s, such as the
Theorems of Synge, Bonnet-Myers, Rauch [1], Berger [2], and Klingenberg [3],
are proved by pertinent, clever applications of Jacobi-Morse theory to
analyze the global geometry of geodesics in the given $M^{n}$ with the kinds
of positivity of curvature assumptions, namely the arc-length element $%
ds^{2} $ remained fixed; while the major works on this topic occurring after
1982 have almost entirely adapted the route of \textit{Ricci evolutions }of
the arc-length elements $ds_{t}^{2}$ introduced in 1982 in Hamilton [4].

Note that the \emph{validity of such proofs} of those major works, of
course, depends entirely on the \emph{preservation }of the kind of \emph{%
curvature positivity assumption }under Ricci evolution of ($M^{n},ds^{2}$),
while the existence of "counterexamples" exhibiting \emph{non-preservation}
\ of the kind of curvature positivity assumption proves the \emph{invalidity
of such proofs.} Anyhow, a \ down-to-earth reality check of how those
absolutely crucial preservation of curvature-positive properties were
proved, as well as the actual validities, must be the\emph{\ first step of
understanding of those works since 1982}.

To begin with, in view of the Ricci evolution equation consisting of a
system of $n(n+1)/2$ of highly nonlinear parabolic equations, one naturally
expects such a proof must be based upon some in-depth, pertinent analysis of
such a system of PDEs of \emph{non-linear diffusion type.} However, it is
disappointing that we failed to find such analytical understanding on the
one hand, and on the other hand, in all of those papers, the indispensable
proofs of \emph{preservation of curvature positivity} are based upon the
kind of "maximum principle" invented by Hamilton [5] in the setting of his 
\emph{reformulated} PDE system instead of the original Ricci evolution
equations. Therefore, such a reality check on those works led us to the
validity (??) of such a reformulation of the Ricci evolution equations.

In short, such a follow-up reality check led us to realize that the
reformulation of Hamilton, in fact, radically changes the basic nature of
the Ricci evolution equations, and hence the proofs of preservation of those
curvature-positivity assumptions under Ricci flows are \emph{invalid} !
However, the above efforts of reality-checking led us to realize that, to 
\emph{avoid misunderstanding of Ricci evolution, }one should begin with a
down-to-earth analysis of simple families of examples, instead of wishful
reformulations without any real analytical understanding beyond the totally
trivial one of $M_{0}^{n}=S^{n}(R)$, and moreover, it is better to be
open-minded to seek the \emph{reality in both directions}, namely either
proving or constructing counterexamples.

Such a viewpoint of geometric analysis naturally led us to, firstly, study
the special case of Ricci evolutions of $O(n)$-invariant hypersurfaces in $%
\mathbb{E}^{n+1}$, namely $C^{\infty }-O(n)$-invariant generalized
rotational hypersurfaces which are generated by $C^{\infty }$ curves in $%
\mathbb{E}_{0}^{2}\subset $ $\mathbb{E}^{n+1},\gamma _{0}\subset \mathbb{E}%
^{2}$, whose analytic geometries can be entirely and neatly characterized by
their \emph{curvature profiles}, that is, the function $k_{1}(\sigma ),a\leq
\sigma \leq b$, recording the curvature of $\gamma _{0}$ at $\gamma
_{0}(\sigma )$, which uniquely determines the curve $\gamma _{0}\subset 
\mathbb{E}^{2}(x,y)$ modulo $x$-translation (for a fixed initial value of $y$%
).

Geometrically, they are straightforward extensions of the classical
rotational surfaces $M^{2}(\gamma _{0})\subset \mathbb{E}^{3}$, and 
\begin{equation}
(\gamma _{0},\mathbb{E}_{0}^{2})\subset (M^{2}(\gamma _{0}),\mathbb{E}%
^{3})\subset (M^{n}(\gamma _{0}),\mathbb{E}^{n+1})  \tag{1.1}
\end{equation}

\noindent clearly shows that the geometric invariant theory of $M^{n}(\gamma
_{0})\subset \mathbb{E}^{n+1}$ can be simply and neatly reduced to that of $%
(M^{2}(\gamma _{0}),\mathbb{E}^{3})$, which are, clearly, totally geodesic
subsurfaces of $M^{n}(\gamma _{0})\subset \mathbb{E}^{n+1}$ (cf. Section 2),
thus providing a neat extension of the Ricci flow equation of $M^{2}(\gamma
_{0})\subset \mathbb{E}^{3})$ to that of $M^{n}(\gamma _{0})\subset \mathbb{E%
}^{n+1}$ (cf. Sections 3 and 4), which can be wholesomely summarized into
the Ricci evolution equation of the curvature profile $k_{1}(\sigma )$ (cf.
Equation (4.14) of Lemma 3 and 3'), a surprising and remarkably simple
geometric invariant theory of $O(n)$-equivariant Ricci evolutions, which, in
turn, naturally led us to another surprising discovery, e.g. the existence
of remarkably ordinary and very simple $M^{n}(\gamma _{0})$ with the
curvature profile just an almost infinitesimal modification of the constant
function of $1/R$, $R\gg 1$ (cf. Section 5), whose Ricci evolution is \emph{%
no longer positive everywhere} for $0\leq t\leq \frac{4}{7}(\frac{1}{R^{2}}+%
\frac{1}{R^{3}}$), namely, the existence of counterexamples toward the
proofs of all those major works since 1982 such as [6, 7, 8, 9, 10, 11, 12,
13] and hence also constituting \emph{counterexamples to the validity of
Hamilton's maximum principle \ for Ricci evolution equations }(cf. Section
6).

\section{Basic differential geometry of rotational surfaces in $\mathbb{E}%
^{3}$ and $O(n)$-hypersurfaces in $\mathbb{E}^{n+1}$}

Let $\mathbb{E}^{n+1}$ be the $(n+1)$-dimensional Euclidean space and $M^{n}$
be an $O(n)$-invariant $C^{\infty }$-hypersurface with $\mathbb{E}_{0}^{1}$
as the axis and consider the inclusions%
\begin{equation*}
\begin{array}{ccccccc}
\mathbb{E}^{n+1} & \supset & \mathbb{E}_{0}^{3} & \supset & \mathbb{E}%
_{0}^{2} & \supset & \mathbb{E}_{0}^{1} \\ 
\cup &  & \cup &  & \cup &  &  \\ 
M^{n} & \supset & \Sigma _{0}^{2} & \supset & \gamma _{0} &  & 
\end{array}%
\end{equation*}%
where $\Sigma _{0}^{2}=M^{n}\cap \mathbb{E}_{0}^{3}$ and $\gamma
_{0}=M^{n}\cap \mathbb{E}_{0}^{2}$ are, respectively, the restricted
rotational surface and the generating curve of both $M^{n}$ and $\Sigma
_{0}^{2}$, and, moreover, $\gamma _{0}$ is a $C^{\infty }$-plane curve in $%
\mathbb{E}_{0}^{2}$ modulo $x$-translations with respect to a Cartesian
coordinate system $(x,y)$ in $\mathbb{E}_{0}^{2}$.

Note that the $O(n)$-symmetry of $M^{n}$ implies, on the one hand, that $%
\Sigma _{0}^{2}$ is the fixed point set of a subisometric group and hence is
a totally geodesic submanifold of $M^{n}$, and, on the other hand, the
choice of $\mathbb{E}_{0}^{3}$ only amounts to the choice of a third one in
the extension of the orthonormal basis of $\mathbb{E}_{0}^{2}$ to that of $%
\mathbb{E}^{n+1}:$

\begin{equation}
\mathbb{E}_{0}^{2}=\langle \mathbf{e}_{0},e_{1}\rangle \subset \mathbb{E}%
_{0}^{3}=\langle \mathbf{e}_{0},e_{1},e_{2}\rangle ...\subset \mathbb{E}%
^{n+1}  \tag{2.1}
\end{equation}%
Therefore, the intrinsic differential geometry of $M^{n}$ is a
straightforward extension of that of $\mathbb{E}_{0}^{2}$. In particular, it
is natural to use the following polar coordinate system for $\mathbb{E}%
_{0}^{2}$ and the generalized spherical coordinate system for $M^{n}$,
namely,

\begin{equation}
\Sigma _{0}^{2} =\{(\sigma ,\theta ):ds^{2}=d\sigma ^{2}+y^{2}d\theta ^{2}\}
\tag{2.2}
\end{equation}

\begin{equation}
M^{n} = \{(\sigma ,\psi _{(n-1)}):ds^{2}=d\sigma ^{2}+y^{2}d\psi
_{(n-1)}^{2}\}  \notag
\end{equation}%
where $\sigma $ is the arc-length parameter of $\gamma _{0}$, $y=y(\sigma )$%
, and $d\psi _{(n-1)}^{2}$ is the spherical metric on $S^{n-1}(1).$

\subsection{Basics on rotational surfaces of $\mathbb{E}_{0}^{3}$}

\begin{figure}[h!]
\centering
\begin{tikzpicture}[scale=0.6, transform shape]
\draw[line width=0.4pt,green] (3.56,5.38) -- ++(-32.4:8.22cm);
\draw[line width=1pt] (7,1) arc (0:120:4);
\draw[black][-Stealth] (5, 4.46) -- (3.85, 2.65);   
\draw[black][-Stealth] (5, 4.44) -- (5, 2.3);
\draw[black, dashed][-Stealth] (5, 4.46) -- (3.5, 5.42);
\draw[->,thick] (1,1)--(13,1) node[right]{\Large${x}$};
\draw[->,thick] (2,0)--(2,8) node[above]{\Large$y$};
\draw[line width=0.4pt, black] (11,1) arc (0:136:0.78cm);
\node at (10.5,2) {$(\pi-\lambda)$};
\node at (5,4.45) {$\bullet$};
\node at (4.2,5.2) {$\mathbf{t}$};
\node at (4.2,3.6) {$\mathbf{n}$};
\node at (5.2,4.7) {$P_0$};
\node at (5.45,3.3) {$\displaystyle{-\frac{\partial}{\partial y}}$};
    \coordinate (A) at (4.5,3.7);  
    \coordinate (B) at (5,3.6);    
     \draw[line width=0.4pt, black] (A) to[bend right=45] (B);
     \node at (4.8,3.8) {$\lambda$};
     \end{tikzpicture}
\caption{The generating curve $\protect\gamma_0$}
\label{fig:RaziPlot1}
\end{figure}
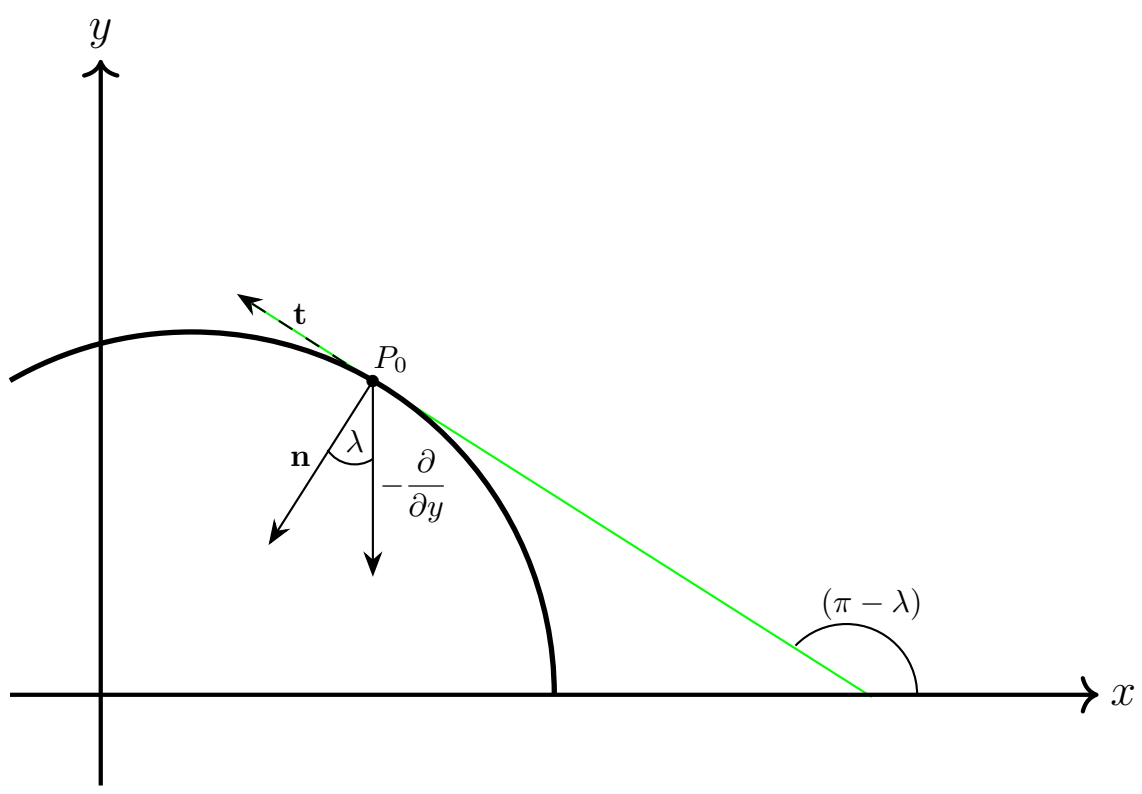

Let $\gamma _{0}\subset $ $\mathbb{E}_{0}^{2}$ be a $C^{\infty }$-curve in $%
\mathbb{E}_{0}^{2}(x,y)$ parametrized by the arc-length $\sigma $ and with
its starting point $P_{0}$ perpendicular to the $x$-axis. As indicated in
Figure~\ref{fig:RaziPlot1} 1, set $(\pi -\lambda )$ as the angle $\angle 
\frac{\partial }{\partial x},\mathbf{t})$ and $\lambda =\angle \mathbf{n,-}%
\frac{\partial }{\partial y})$, namely,

\begin{equation}
\begin{array}{c}
\gamma _{0}=\left\{ \gamma _{0}(\sigma )=(x(\sigma ),y(\sigma )),0\leq
\sigma \leq l\right\} \\ 
\mathbf{t(}\sigma )=(\frac{\partial x}{d\sigma },\frac{\partial y}{d\sigma }%
)=(-\cos \lambda (\sigma ),\sin \lambda (\sigma )) \\ 
\mathbf{n(}\sigma )=(-\frac{\partial y}{d\sigma },\frac{\partial x}{d\sigma }%
)=(-\sin \lambda (\sigma ),-\cos \lambda (\sigma ))%
\end{array}
\tag{2.3}
\end{equation}

\bigskip

\bigskip

Therefore, the curvature vector at $\gamma _{0}(\sigma )$ of the $y(\sigma )$%
-circle is given by $(0,-\frac{1}{y})$ and hence, the principal curvatures
of $\gamma _{0}(\sigma )$ and $M^{2}(\gamma _{0})$ are given by Meusnier's
formula, namely 
\begin{equation}
k_{1}=-\frac{\partial \lambda }{\partial \sigma }\text{ \ and }k_{2}=\frac{%
\cos \lambda (\sigma )}{y}=-\frac{1}{y}\frac{\partial x}{\partial \sigma } 
\tag{2.4}
\end{equation}%
with $\mathbf{e}_{1}=\frac{\partial }{\partial \sigma }$ and $\mathbf{e}_{2}=%
\frac{1}{y}\frac{\partial }{\partial \theta }$ as their respective principal
directions.



Note that%
\begin{equation}
(\frac{\partial \mathbf{t}}{\partial \sigma })\cdot \mathbf{t} =\frac{%
\partial ^{2}x}{\partial \sigma ^{2}}\frac{\partial x}{\partial \sigma }+%
\frac{\partial ^{2}y}{\partial \sigma ^{2}}\frac{\partial y}{\partial \sigma 
}=0 \notag
\end{equation}

\begin{equation}
k_{1} =\frac{\partial \mathbf{t}}{\partial \sigma }\cdot \mathbf{n}=-\frac{%
\partial ^{2}x}{\partial \sigma ^{2}}\frac{\partial y}{\partial \sigma }+%
\frac{\partial ^{2}y}{\partial \sigma ^{2}}\frac{\partial x}{\partial \sigma} \notag
\end{equation}%

Hence, 
\begin{equation}
k =k_{1}k_{2}=k_{1}\frac{\cos \lambda (\sigma )}{y}=\frac{1}{y}\frac{%
\partial x}{\partial \sigma }\left\{ \frac{\partial ^{2}x}{\partial \sigma
^{2}}\frac{\partial y}{\partial \sigma }-\frac{\partial ^{2}y}{\partial
\sigma ^{2}}\frac{\partial x}{\partial \sigma }\right\} \notag \\
\end{equation}
\begin{equation}
=-\frac{1}{y}\frac{\partial ^{2}y}{\partial \sigma ^{2}}\left\{ (\frac{%
\partial y}{\partial \sigma })^{2}+(\frac{\partial x}{\partial \sigma }%
)^{2}\right\} =-\frac{1}{y}\frac{\partial ^{2}y}{\partial \sigma ^{2}} \tag{2.5}
\end{equation}

\bigskip\ \ \ \ \ 

\ \ 

\noindent \textbf{Remarks 1}

\textbf{(i)} The above formula of intrinsic curvature of $M_{0}^{2}$ is the
specialization of an abstract rotational surface, that is, $(O(2),M^{2}(f)):$%
\begin{equation}
M^{2}(f)=\left\{ (\sigma ,\theta ):ds^{2}=d\sigma ^{2}+f(\sigma )^{2}d\theta
^{2}\right\}  \tag{2.6}
\end{equation}%
where $\sigma $ is the arc-length parameter of the orbit space $%
M^{2}(f)/O(2) $, while $f(\sigma )$ is the size function, namely the length
of the $O(2)$-orbit $G(\sigma )$ is given by $2\pi f(\sigma )$.

\textbf{(ii)} Set $\alpha =\angle \frac{\partial }{\partial \sigma },\mathbf{%
t}).$ Then $\frac{\partial }{\partial \theta }$ is the Killing vector field
and $f(\sigma )=\left\vert \frac{\partial }{\partial \theta }\right\vert $,
while the \emph{generalized sine law }(i.e. the conservation law or
Noether's law) is that%
\begin{equation}
\left\langle \frac{\partial }{\partial \theta },\mathbf{t}\right\rangle
=f(\sigma )\sin \alpha =\text{ constant,}  \tag{2.7}
\end{equation}%
and it follows from the Gauss-Bonnet formula that the Gaussian curvature of $%
M^{2}(f)$ is give by

\begin{equation}
k=-\frac{f^{\prime \prime }}{f}  \tag{2.8}
\end{equation}

All of these are basic consequences of the first variational formula of
arc-length which implies the following formula of \emph{geodesic curvature: }

\begin{equation}
k_{g}(\gamma (s))=\frac{\partial \alpha }{\partial s}+f^{\prime }(\sigma )%
\frac{\partial \sigma }{\partial s}  \tag{2.9}
\end{equation}

\textbf{(iii)} It is interesting to note here that in the special case of $%
ISO(O(2),M^{2}(f))\simeq O(2)$ and $\pi _{1}(M^{2}(f))=1$, one has

\begin{equation*}
k=\text{const. }(k_{0})\Longleftrightarrow \lbrack f^{\prime \prime
}+k_{0}f=0,f(0)=0,f^{\prime }(0)=1]
\end{equation*}%
and this is equivalent to%
\begin{equation}
f(\sigma )=\ 
\begin{array}{cc}
\sigma \text{ } & k_{0}=0 \\ 
\frac{1}{\sqrt{k_{0}}}\sin \sqrt{k_{0}}\sigma & k_{0}>0 \\ 
\frac{1}{\sqrt{-k_{0}}}\sinh \sqrt{-k_{0}}\sigma & k_{0}<0%
\end{array}
\tag{2.10}
\end{equation}%
which provides a unifying theory of Euclidean, spherical, and non-Euclidean
geometry.

\subsection{Basic differential geometry of $O(3)$-hypersurfaces in $\mathbb{E%
}^{4}$}

Note that $M^{3}$ has a family of totally geodesic submanifolds, $%
\{M_{\varphi }^{2}\}$, with $M_{0}^{2}=M^{3}\cap \mathbb{E}_{0}^{3}.$
Therefore, the principal curvatures (directions) of $M^{3}$ are also given
by the above $k_{1}$ and $k_{2}$ at $\gamma _{0}(\sigma )$, while their
corresponding directions are the tangential vector $\frac{\partial }{%
\partial \sigma }$ and $\{\frac{1}{y}\frac{\partial }{\partial \theta }$
parametrized by $\varphi \}$. Hence, it is natural to choose the orthonormal
frame field at $\gamma _{0}(\sigma )$ as follows, namely%
\begin{equation}
\mathbf{e}_{1}=\frac{\partial }{\partial \sigma }\text{, \ \ }e_{2}=\frac{1}{%
y}\frac{\partial }{\partial \theta }|_{\varphi =0}\text{, \ \ }e_{3}=\frac{1%
}{y}\frac{\partial }{\partial \theta }|_{\varphi =\frac{\pi }{2}}  \tag{2.11}
\end{equation}%
while%
\begin{equation}
X_{2}=\frac{\partial }{\partial \theta }|_{\varphi =0}\text{, \ }X_{2}=\frac{%
\partial }{\partial \theta }|_{\varphi =\frac{\pi }{2}}  \tag{2.12}
\end{equation}%
are an orthogonal pair of KIlling vector fields on $M^{3}$. Moreover, the
Ricci tensor of $M^{3}$ and the principal curvatures are related as follows,
namely

\begin{equation}
Ric(\frac{\partial }{\partial \sigma },\frac{\partial }{\partial \sigma }) =
2k_{1}k_{2}=-\frac{2}{y}\frac{\partial ^{2}y}{d\sigma ^{2}}  \notag \\
\end{equation}
\begin{equation}
Ric(\mathbf{e}_{2,}\mathbf{e}_{2})=Ric(\mathbf{e}_{3,}\mathbf{e}_{3}) =
(k_{1}+k_{2})k_{2}=\frac{1}{y^{2}}\{1-\frac{d}{d\sigma }(y.\frac{dy}{d\sigma 
})\}  \tag{2.13} \\
\end{equation}
\begin{equation}
Ric(X_{2},X_{2}) = Ric(X_{3},X_{3})= y^{2}Ric(\mathbf{e}_{2,}\mathbf{e}%
_{2})= \{1-\frac{d}{d\sigma }(y.\frac{dy}{d\sigma })\}  \notag
\end{equation}

\noindent \textbf{Remark 2} Straightforward generalization will show that
for $M^{n},n>3$, one has 
\begin{equation}
Ric(\frac{\partial }{\partial \sigma },\frac{\partial }{\partial \sigma }%
)=(n-1)k_{1}k_{2} =-\frac{(n-1)}{y}\frac{\partial ^{2}y}{d\sigma ^{2}} 
\tag{2.14a} \\
\end{equation}
\begin{equation}
Ric(\mathbf{e}_{j},e_{j})=k_{1}k_{2}+(n-2)k_{2}^{2}=-\frac{1}{y}\frac{%
\partial ^{2}y}{d\sigma ^{2}}+(n-2)\frac{\cos ^{2}\lambda }{y^{2}} 
\tag{2.14b}
\end{equation}

\section{Intrinsic characterization of $O(n)$-hypersurfaces in $\mathbb{E}%
^{n+1}$}

In the general setting of Riemannian manifolds, an $O(n)$-manifold $M^{n}$
can also be represented by its spherical coordinate system, namely%
\begin{equation}
M^{n}=\{(\sigma ,\psi _{(n-1)}):ds^{2}=d\sigma ^{2}+f(\sigma )^{2}d\psi
_{(n-1)}^{2}  \tag{3.1}
\end{equation}%
where $\sigma $ is the arc-length parameter of the 1-dimensional orbit space 
$M^{n}/O(n)$ and $f(\sigma )$ is the size function $G(\sigma )\simeq
S^{n-1}(f$($\sigma )).$

In the special case of $O(n)$-hypersurfaces in $\mathbb{E}^{n+1}$, $%
M^{n}(\gamma )/O(n)\simeq \gamma $, $\gamma $ is the generating curve in $%
\mathbb{E}^{2}(x,y)$, $f(\sigma )=y(\sigma )$. Therefore, one of the \emph{%
necessary conditions } that $M^{n}(f(\sigma ))$ can be realized as an $O(n)$%
-hypersurface is that 
\begin{equation}
\left( \frac{d}{d\sigma }f(\sigma )\right) ^{2}\leq (\frac{dy}{d\sigma }%
)^{2}+(\frac{dx}{d\sigma })^{2}=1  \tag{3.2}
\end{equation}%
and moreover, it is also quite obvious that $f(\sigma )$ can \emph{not} have
maximal (resp. minimal) points with $\frac{df}{d\sigma }|_{\hat{\sigma}_{\mp
}\text{ \ }}$equal to $\pm 1$ (resp. $\mp 1$ with nonzero $f(\hat{\sigma})$).

\begin{lemma}
An $O(n)$-manifold $M$ with size function $f(\sigma )$ can be realized as an 
$O(n)$-hypersurface if and only if 
\begin{equation}
\left( \frac{df}{d\sigma }\right) ^{2}\leq 1  \tag{3.3}
\end{equation}%
and $f(\sigma )$ has no maximal (resp. non-zero minimal) points with 
\begin{equation}
\frac{df}{d\sigma }|_{\hat{\sigma}_{\mp }\text{ \ }}=\pm 1\text{ }(resp\text{%
. }\mp 1)  \tag{3.4}
\end{equation}
\end{lemma}

\begin{proof}
The above condition is clearly necessary. Thus, it suffices to prove the
sufficiency by construction.$\ $

Note that the smoothness of $M^{n}$ implies that, in the vicinity of the
fixed points of $M^{n}$, which would be the terminating (resp. starting )
points on the $x$-axis for its realization, we must have $\frac{df}{d\sigma }%
=+1$ or $-1$.

Let us first consider the case that there is no point of $\sigma $ with $%
\frac{df}{d\sigma }=$ $\pm 1$. Then one simply constructs the generating
curve $\gamma (\sigma )$ by setting 
\begin{equation}
y(\sigma )=f(\sigma )\text{ \ \ and \ }x(\sigma )=\int_{\sigma _{0}}^{\sigma
}\pm \sqrt{1-(\frac{df}{d\sigma })^{2}}d\sigma  \tag{3.5}
\end{equation}%
where the choice of sign will give two generating curves reflectionally
symmetric to each other.

In general, the whole range of $\sigma $ can be divided into segments of the
above kind by points or close intervals with $\frac{df}{d\sigma }$
identically equal to $\pm 1$, whose corresponding $\gamma (\sigma )$ will be
just vertical intervals. Therefore, it is not difficult to see that the
entire generating curve can be obtained just by the piecing together of such
segments.
\end{proof}

\noindent \textbf{Remark 3 }Note that such a $C^{\infty }$ plane curve $%
\gamma _{t}$ modulo $x$-translations is uniquely characterized \ by its
curvature profile together with its minimum value of $y$, say $y_{0}$. Thus, 
$M^{n}$ will be diffeomorphic to $S^{n-1}\times I$ if $y_{0}\neq 0$ or $%
\frac{\partial y}{\partial \sigma }=\pm 1$ if $y_{0}=0$ (i.e. at its fixed
points). Therefore, in the viewpoint (or spirit) of \emph{intrinsic
geometric invariant theory, }the corresponding \emph{evolution equation of
the curvature profiles, }$k_{1}(\gamma _{t}(\sigma ))$, should be the
fundamental key for the study of such Ricci flows (cf. Section 4).\emph{\ }

\section{Ricci flows of $O(n)$-hypersurfaces and the equation of Ricci
evolution of their generating curves}

\subsection{The Ricci flow and equations of $O(3)$-hypersurfaces in $\mathbb{%
E}^{4}$.}

Let $\{$ $M_{t}^{3},0\leq t\leq T\}$ be a Ricci flow of $O(3)$-hypersurfaces
with $\{\gamma _{t},0\leq t\leq T\}$ as their generating curves, namely%
\begin{equation}
M_{t}^{3}=\{(u,\theta ,\varphi):
ds_{t}^{2}=g(t,u)^{2}du^{2}+y(t,u)^{2}(d\theta ^{2}+\sin ^{2}\theta d\varphi
^{2})\}  \tag{4.1}
\end{equation}%
where $(g(t,u),y(t,u))$ satisfy the following evolution equation:%
\begin{equation}
(\mathcal{E}_{0}):\text{ \ \ \ \ \ }\left\{ 
\begin{array}{c}
\frac{\partial }{\partial t}[g(t,u)^{2}]=-2Ric^{t}(\frac{\partial }{\partial
u},\frac{\partial }{\partial u}) \\ 
\frac{\partial }{\partial t}[y(t,u)^{2}]=-2Ric^{t}(\frac{\partial }{\partial
\theta },\frac{\partial }{\partial \theta })%
\end{array}%
\right.  \tag{4.2}
\end{equation}%
and where $Ric^{t}(\cdot ,\cdot )$ is the Ricci tensor of $M_{t}^{3}$ along $%
\gamma _{t}$ (i.e. at $\theta =0$) and the arc-length parameter $\sigma _{t} 
$ of $\gamma _{t}$ is given by%
\begin{equation}
d\sigma _{t}=g(t,u)du\text{ \ \ or \ }\frac{\partial }{du}=g(t,u)\frac{%
\partial }{d\sigma _{t}}  \tag{4.3}
\end{equation}

Now, let us make use of the basic geometry of Section 2 to compute the
reduction of $(\mathcal{E}_{0})$ to that of the evolution equation of Ricci
evolution of $\{\gamma _{t}\}$. Set%
\begin{equation}
\rho (\sigma _{t})=y(t,u)\text{ \ and }d\sigma _{t}=g(t,u)du\text{ (as in
(4.3))}  \tag{4.4}
\end{equation}%
Then%
\begin{equation}
ds_{t}^{2} =d\sigma _{t}^{2}+\rho (\sigma _{t})^{2}(d\theta ^{2}+\sin
^{2}\theta d\varphi ^{2})  \tag{4.5} \\
\end{equation}
\begin{equation}
\frac{\partial }{\partial u} = \frac{d\sigma _{t}}{du}\frac{\partial }{%
\partial \sigma _{t}}=g(t,u)\frac{\partial }{\partial \sigma _{t}}\text{ \
(i.e. chain rule)}  \notag
\end{equation}

\noindent where $g(0,u)\equiv 1$, that is, $u$ is the arc-length of $\gamma
_{0}$ to begin with. By the formula of Section 2, adapting to $\gamma _{t}$
and $\sigma _{t}$, one has 
\begin{equation}
\begin{array}{c}
Ric(\frac{\partial }{\partial u},\frac{\partial }{\partial u})=g^{2}Ric^{t}(%
\frac{\partial }{\partial \sigma _{t}},\frac{\partial }{\partial \sigma _{t}}%
)=g^{2}(-2\frac{\rho ^{\prime \prime }}{\rho }) \\ 
Ric(\frac{\partial }{\partial \theta },\frac{\partial }{\partial \theta }%
)=\rho ^{2}Ric^{t}(\mathbf{e}_{2},\mathbf{e}_{2})=1-(\rho \cdot \rho
^{\prime })^{\prime }%
\end{array}
\tag{4.6}
\end{equation}

Set 
\begin{equation}
G=g(t,u)^{2},\ Y=y(t,u)^{2}  \tag{4.7}
\end{equation}%
It is advantageous to rewrite the evolution equation as follows, namely

\begin{lemma}
In the special case of Ricci flow \ of $O(3)$-hypersurfaces, the system $(%
\mathcal{E}_{0})$ can be rewritten neatly in terms of the geometric
invariants of the generating curves $\{\gamma _{t}\}$, namely%
\begin{equation}
(\mathcal{E}_{1}):\left\{ 
\begin{array}{c}
\frac{\partial }{\partial t}\ln G=-4k_{1}k_{2}\text{, \ }G=g(t,u)^{2} \\ 
\frac{\partial }{\partial t}\ln Y=-2(k_{1}+k_{2})k_{2}\text{, \ }Y=y(t,u)^{2}%
\end{array}%
\right\}  \tag{4.8}
\end{equation}%
where $k_{1}$ is the curvature of the plane curve $\gamma _{t}$, and $%
k_{2}=\cos \lambda _{t}/y=-\frac{1}{y}\frac{dx}{d\sigma _{t}}$ \ (c.f.
(2.4)).
\end{lemma}

\begin{proof}
By direct substitution and $(\mathcal{E}_{0})$, note that%
\begin{equation}
\frac{\partial }{\partial t}\ln G =\frac{\partial }{\partial t}G/G=-\frac{2}{%
G}Ric^{t}(\frac{\partial }{\partial u},\frac{\partial }{\partial u}) 
\tag{4.9a} \\
\end{equation}

\begin{equation}
= -2Ric^{t}(\frac{\partial }{\partial \sigma _{t}},\frac{\partial }{\partial
\sigma _{t}})=-4k_{1}k_{2}  \notag
\end{equation}

\begin{equation}
\frac{\partial }{\partial t}\ln Y = \frac{\partial }{\partial t}Y/Y=-\frac{2%
}{Y}Ric^{t}(\frac{\partial }{\partial \theta },\frac{\partial }{\partial
\theta })  \tag{4.9b} \\
\end{equation}
\begin{equation}
= -2Ric^{t}(\mathbf{e}_{2},\mathbf{e}_{2})=-2(k_{1}+k_{2})k_{2}  \notag
\end{equation}
\end{proof}

\textbf{Corollary 1}

\begin{equation}
(\mathcal{E}_{2}):\text{ \ \ \ \ \ }\left\{ 
\begin{array}{c}
\frac{\partial }{\partial t}g=-2k_{1}k_{2}g\text{ \ \ \ \ \ \ \ \ \ \ \ \ \
\ \ \ \ \ \ \ \ \ \ \ \ \ \ \ \ \ \ \ \ \ \ \ \ \ \ \ \ \ \ \ \ \ \ \ \ \ \
\ \ \ \ \ \ } \\ 
\frac{\partial }{\partial t}y=-yk_{2}(k_{1}+k_{2})=-k_{1}\cos \lambda -\frac{%
\cos ^{2}\lambda }{y}=(\frac{\partial }{\partial \sigma _{t}})^{2}y-\frac{%
\cos ^{2}\lambda }{y}%
\end{array}%
\right.  \tag{4.10}
\end{equation}

\textbf{Corollary 2 } 
\begin{equation}
(\mathcal{E}_{2})_{Y}:\left\{ 
\begin{array}{l}
\text{\ }\frac{\partial }{\partial t}Y=-2Y\text{ }Ric^{t}(\mathbf{e}_{2},%
\mathbf{e}_{2})=-2Ric^{t}(\frac{\partial }{\partial \theta },\frac{\partial 
}{\partial \theta })\text{\ \ \ \ \ \ \ \ \ \ \ \ \ \ \ \ \ \ \ \ \ \ \ \ \
\ \ \ \ \ \ \ \ \ \ \ \ \ \ \ \ \ \ \ \ \ \ \ \ \ \ \ \ \ \ \ \ \ \ \ \ } \\ 
\text{ \ \ \ \ \ \ }=-2(1-(\rho \cdot \rho ^{\prime })^{\prime }=(\rho
^{2})^{\prime \prime }-2\text{ \ \ \ \ } \\ 
\text{ \ \ \ \ \ \ }=(\frac{\partial }{\partial \sigma _{t}})^{2}Y-2\text{\
\ \ \ \ \ \ \ \ \ \ \ \ \ \ \ \ \ \ \ \ \ \ \ \ \ \ \ \ \ \ \ \ \ \ \ \ \ \
\ \ \ \ \ \ \ \ \ \ \ \ \ \ \ \ \ }%
\end{array}%
\right.  \tag{4.11}
\end{equation}

\noindent \textbf{Remarks 4}

\textbf{(i)} $(\mathcal{E}_{2})$ and $(\mathcal{E}_{2})_{Y}$ are remarkably
simple non-linear parabolic equations of space dimension 1, which are simple
natural generalizations of the \emph{linear} heat equation. Physically
speaking, it is the diffusion equation along a wire whose "heat
conductivity" is no longer constant, but is rather time dependent and given
by the equation of $\frac{\partial }{\partial t}g$.

\textbf{(ii)} The \emph{non-linearity of }$(\mathcal{E}_{2})_{Y}$ solely
lies in the time dependence of $\sigma _{t}$, while that of $(\mathcal{E}%
_{2})$ has another non-linear term "$-\cos ^{2}\lambda /y$ \ ".

\textbf{(iii)} As mentioned in Section 1, it is conceptually intrinsic and
technically advantageous to have the \emph{evolution equation of the
curvature profiles }$k_{1}^{t}$ of the generating curves $\gamma _{t}$,
which is actually a corollary of Lemma 2.

\bigskip Now, note first that%
\begin{equation*}
\frac{\partial \lambda }{\partial \sigma _{t}}=-k_{1}^{t}\text{ (see }(2.4)%
\text{)}
\end{equation*}%
\ 
\begin{equation}
\sin \lambda \frac{\partial \lambda }{\partial t}k_{1}^{t}-\cos \lambda 
\frac{\partial k_{1}^{t}}{\partial t} = \frac{\partial }{\partial t}\cdot 
\frac{\partial ^{2}}{\partial \sigma _{t}^{2}}y=\frac{\partial ^{2}}{%
\partial \sigma _{t}^{2}}\frac{\partial y}{\partial t} \   \tag{4.12a}
\end{equation}
\begin{equation}
= \frac{\partial ^{2}}{\partial \sigma _{t}^{2}}(-\cos \lambda k_{1}^{t})-%
\frac{\partial ^{2}}{\partial \sigma _{t}^{2}}(\frac{\cos ^{2}\lambda }{y}) 
\notag
\end{equation}
\begin{equation}
= \left[ \cos \lambda (k_{1}^{t})^{2}-\sin \lambda \frac{\partial k_{1}^{t}}{%
\partial \sigma _{t}}\right] k_{1}^{t}-2\sin \lambda k_{1}^{t}\cdot \frac{%
\partial k_{1}^{t}}{\partial \sigma _{t}}  \notag \\
-\cos \lambda \frac{\partial ^{2}k_{1}^{t}}{\partial \sigma _{t}^{2}}-\frac{%
\partial ^{2}}{\partial \sigma _{t}^{2}}(\frac{\cos ^{2}\lambda }{y})  \notag
\end{equation}
\newline

\noindent hence, by solving for $\frac{\partial k_{1}^{t}}{\partial t}$ from
Eq. 4.12a one has \newline
\begin{equation}
\frac{\partial k_{1}^{t}}{\partial t}=\tan \lambda \frac{\partial \lambda }{%
\partial t}k_{1}^{t}-(k_{1}^{t})^{3}+3\tan \lambda k_{1}^{t}\frac{\partial
k_{1}^{t}}{\partial \sigma _{t}}+\frac{\partial ^{2}k_{1}^{t}}{\partial
\sigma _{t}^{2}}+\sec \lambda \frac{\partial ^{2}}{\partial \sigma _{t}^{2}}(%
\frac{\cos ^{2}\lambda }{y})  \tag{4.12b}
\end{equation}
\newline

Next, let us compute $\tan \lambda \frac{\partial \lambda }{\partial t}$ via
differentiation of $\sin \lambda =\frac{\partial y}{\partial \sigma _{t}}$;
one has%
\begin{equation}
\cos \lambda \frac{\partial \lambda }{\partial t} = \frac{\partial ^{2}y}{%
\partial t\partial \sigma _{t}}=\frac{\partial }{\partial \sigma _{t}}%
\left\{ -\cos \lambda k_{1}^{t}-\frac{\cos ^{2}\lambda }{y}\right\} 
\tag{4.13a}
\end{equation}

\begin{equation}
= -\sin \lambda (k_{1}^{t})^{2}-\cos \lambda \frac{\partial k_{1}^{t}}{%
\partial \sigma _{t}}-\frac{\partial }{\partial \sigma _{t}}(\frac{\cos
^{2}\lambda }{y})  \notag
\end{equation}
\newline

\noindent thus having 
\begin{equation}
\tan \lambda \frac{\partial \lambda }{\partial t}=-\tan ^{2}\lambda
(k_{1}^{t})^{2}-\tan \lambda \frac{\partial k_{1}^{t}}{\partial \sigma _{t}}%
-\sec \lambda \tan \lambda \frac{\partial }{\partial \sigma _{t}}(\frac{\cos
^{2}\lambda }{y})  \tag{4.13b}
\end{equation}%
Therefore, a simple substitution of the RHS of Eq. 4.13b into the first term
of the RHS of Eq. 4.12b yields

\begin{lemma}
The evolution equation of $k_{1}^{t}$ of the generating curves $\gamma _{t}$
of Ricci flow $M^{3}(\gamma _{t})$ for short times is as follows:%
\begin{equation}
(\mathcal{E}_{3}):\frac{\partial k_{1}^{t}}{\partial t}=\frac{\partial ^{2}}{%
\partial \sigma _{t}^{2}}k_{1}^{t}+2\tan \lambda k_{1}^{t}\frac{\partial
k_{1}^{t}}{\partial \sigma _{t}}-\sec ^{2}\lambda (k_{1}^{t})^{3}  \tag{4.14}
\end{equation}
\begin{equation}
+ \sec \lambda \left\{ \frac{\partial ^{2}}{\partial \sigma _{t}^{2}}(\frac{%
\cos ^{2}\lambda }{y})-\tan \lambda \text{ }k_{1}^{t}\frac{\partial }{%
\partial \sigma _{t}}(\frac{\cos ^{2}\lambda }{y})\right\}  \notag
\end{equation}
\end{lemma}

\noindent \textbf{Remarks 5} \textbf{(i)} In the general setting of Ricci
flows of closed 3-manifolds, it is more natural to formulate our study of 
\emph{Ricci flows of }$O(3)$\emph{-manifolds intrinsically} with the
viewpoint of \emph{geometric invariant theory}, while the discussions in
Sections 2, 3, 4 are simple basics of geometric representation theory of a
rather special kind.

\textbf{(ii)} For example, in the formulation of $O(3)$-equivariant Ricci
flows of $3$-manifolds, one has:

\begin{equation*}
\begin{array}{cccc}
(O(3),M_{t}^{3}) & \subset & (O(3),\mathbb{E}^{4}) &  \\ 
\downarrow &  & \downarrow &  \\ 
M^{3}/O(3) & \subset & \mathbb{E}^{4}/O(3)\simeq \mathbb{E}_{+}^{2} & \text{%
(with orbital distance metric)}%
\end{array}%
\end{equation*}%
thus, the generating curves $\gamma _{t}$ are intrinsically the orbit
spaces, $M_{t}^{3}/O(3)$, equipped with the orbital distance metric inside $%
\mathbb{E}^{4}/O(3)\simeq \mathbb{E}_{+}^{2}$, modulo $x$-translations. In
such an intrinsic setting,%
\begin{equation}
\cos \lambda (\sigma _{t})=-\frac{\partial x}{\partial \sigma _{t}}\text{, \
\ }\sin \lambda (\sigma _{t})=\frac{\partial y}{\partial \sigma _{t}},\text{
\ \ }\tan \lambda (\sigma _{t})=-\frac{\partial y}{\partial x}(\sigma _{t}) 
\tag{4.15}
\end{equation}

\noindent and $k_{1}^{t}=-\frac{\partial \lambda }{\partial \sigma _{t}}$,
while the discussion of Section 3 can be extended directly to show that $%
k_{1}^{t}(\sigma _{t})$ also constitutes a complete geometric invariant for $%
M_{t}^{3}/O(3)$ as a curve in $\mathbb{E}_{+}^{2}$ modulo $x$-translations.

\textbf{(iii)} For example, the above setting can also be extended to the
study of $O(3)$-hypersurfaces and their Ricci flows in $4$-dimensional
spherical as well as hyperbolic spaces, etc.

\subsection{The general case of Ricci evolutions of $O(n)$-hypersurfaces in $%
\mathbb{E}^{n+1}$}

In the general case of $n\geq 3$, we have the following

\begin{equation}
(\mathcal{E}_{2}^{\prime }):\left\{ 
\begin{array}{c}
\frac{\partial }{\partial t}g=-(n-1)k_{1}k_{2}g\text{ \ \ \ \ \ \ \ \ \ \ \
\ \ \ \ \ \ \ \ \ \ \ \ \ \ \ \ \ \ \ \ \ \ \ \ \ \ \ \ \ \ \ \ \ \ \ \ \ }
\\ 
\frac{\partial }{\partial t}y=-k_{1}\cos \lambda -(n-2)\frac{\cos
^{2}\lambda }{y}=(\frac{\partial }{\partial \sigma })^{2}y-(n-2)\frac{\cos
^{2}\lambda }{y}%
\end{array}%
\right.  \tag{4.16}
\end{equation}%
thus having (cf. (4.12a))

\begin{equation}
\sin \lambda \frac{\partial \lambda }{\partial t}k_{1}^{t}-\cos \lambda 
\frac{\partial k_{1}^{t}}{\partial t} =\frac{\partial }{\partial t}\cdot 
\frac{\partial ^{2}}{\partial \sigma _{t}^{2}}y  \tag{4.17}
\end{equation}
\begin{equation}
=\frac{\partial ^{2}}{\partial \sigma _{t}^{2}}\frac{\partial y}{\partial t}=%
\frac{\partial ^{2}}{\partial \sigma _{t}^{2}}\left( -k_{1}^{t}\cos \lambda
-(n-2)\frac{\cos ^{2}\lambda }{y}\right)  \notag
\end{equation}
\begin{equation}
= \cos \lambda (k_{1}^{t})^{3}-3\sin \lambda \frac{\partial k_{1}^{t}}{%
\partial t}k_{1}^{t}\ -\frac{\cos ^{2}\lambda }{y}\frac{\partial
^{2}k_{1}^{t}}{\partial \sigma _{t}^{2}}-(n-2)\frac{\partial ^{2}}{\partial
\sigma _{t}^{2}}(\frac{\cos ^{2}\lambda }{y})  \notag
\end{equation}

\noindent and 
\begin{equation}
\cos \lambda \frac{\partial \lambda }{\partial t} = \frac{\partial ^{2}y}{%
\partial t\partial \sigma _{t}}=\frac{\partial }{\partial \sigma _{t}}%
\left\{ -\cos \lambda \cdot k_{1}^{t}-(n-2)\frac{\cos ^{2}\lambda }{y}%
\right\}  \tag{4.18}
\end{equation}
\begin{equation}
= -\sin \lambda (k_{1}^{t})^{2}-\cos \lambda \frac{\partial k_{1}^{t}}{%
\partial \sigma _{t}}-(n-2)\ \frac{\partial }{\partial \sigma _{t}}(\frac{%
\cos ^{2}\lambda }{y})  \notag
\end{equation}

\begin{equation*}
\Longleftrightarrow \text{ }\tan \lambda \frac{\partial \lambda }{\partial t}%
=-\tan ^{2}\lambda (k_{1}^{t})^{2}-\tan \lambda \frac{\partial k_{1}^{t}}{%
\partial \sigma _{t}}-(n-2)\sec \lambda \tan \lambda \frac{\cos ^{2}\lambda 
}{y}
\end{equation*}

Therefore, solving for $\frac{\partial k_{1}^{t}}{\partial t}$ from the
above pair (4.17), (4.18) of equations one has the general form of Lemma 3,
as follows: \newline

\noindent \textbf{Lemma 3'}

\emph{The evolution equation of $k_{1}^{t}$ of the generating curves $\gamma
_{t}$ of $M^{n}(\gamma _{t})\subset \mathbb{E}^{n+1}$ for short time is as
follows:} 
\begin{equation}
(\mathcal{E}_{3}^{\prime }): \frac{\partial k_{1}^{t}}{\partial t}=\frac{%
\partial ^{2}}{\partial \sigma _{t}^{2}}k_{1}^{t}+2\tan \lambda k_{1}^{t}%
\frac{\partial k_{1}^{t}}{\partial \sigma _{t}}-\sec ^{2}\lambda
(k_{1}^{t})^{3}  \tag{4.19}
\end{equation}
\begin{equation}
+(n-2)\sec \lambda \left\{ \frac{\partial ^{2}}{\partial \sigma _{t}^{2}}(%
\frac{\cos ^{2}\lambda }{y})-\tan \lambda \text{ }k_{1}^{t}\frac{\partial }{%
\partial \sigma _{t}}(\frac{\cos ^{2}\lambda }{y})\right\}  \notag
\end{equation}

\noindent \textbf{Remark 6}\newline
In summary, Sections 2, 3, 4 amount to a concise basic theory of Ricci
evolutions of $O(n)$-hypersurfaces in $\mathbb{E}^{n+1},n\geq 3$, while
Lemmas 2 and 3 constitute the fundamental Lemmas. Moreover, such a simple
theory can be readily deduced as an illustrating example of the classical
basic local differential geometry, namely Meusnier's Theorem, Euler's
Theorem, and Gaussian curvature, which amounts to a simple computation of
(2.5), namely%
\begin{equation*}
k=k_{1}k_{2}=-\frac{1}{y}\frac{\partial ^{2}y}{\partial \sigma ^{2}}\text{ \
(cf. (2.5))}
\end{equation*}%
where $k_{1},k_{2}$ are the principal curvatures.\newline

\noindent \textbf{Remark 7}\newline
In retrospect, Ricci tensor as well as Ricci evolutions are conceptually
defined in terms of local coordinate systems; hence it is generally
difficult to provide a global analytical representation, and this is exactly
the root of non-linear difficulties.

However, in the \emph{special} case of $O(n)$-hypersurfaces, the magic
finger of symmetry, in one stroke, makes the local principal directions
already \emph{global}!

\section{Ricci evolutions of a simple family of rather ordinary, small
modifications of spherical $O(3)$ (resp. $O(n)$) hypersurfaces in \{$O(3),%
\mathbb{E}^{4}$\} (resp. \{$O(n),\mathbb{E}^{n+1}$\})}

\subsection{Ricci evolutions -- a misleading venture?}

In this section, as has already been mentioned in Section 1, we shall
proceed to construct a family of simple examples of $O(n)$-equivariant Ricci
evolutions \ for extremely short times in the setting of ($O(n),\mathbb{E}%
^{n+1}$), $n\geq 3$, whose initial $M^{n}(\gamma _{0})$ satisfies the
PRE-property (i.e. positive definite Ricci curvature everywhere), but in an
extremely short, almost infinitesimal time, no longer has the PRE-property!

Therefore, one the one hand, \ they provide simple, easy to check,
counterexamples to the fundamental \emph{PRE-preservation }theorem of R.
Hamilton (1982), [4], as well as the \emph{proclaimed}\ \emph{maximal
principle }[5]\ for Ricci evolution equations that the deduction of the
PRE-preservation theorem is based upon; and on the other hand, such simple
examples and further modifications of them, in fact, indicate that the
entire approach of proving the "Differential Sphere Theorem" via \emph{Ricci
evolutions}\ is \emph{definitely unworkable, }and is a totally misleading
enterprise.\ \ 

\begin{figure}[h!]
\centering
\hspace{-4em} 
\begin{tikzpicture}
\draw[Stealth-] (7.2,1.56) -- (8,2.1);
\draw[black,line width=0.4pt][-] (5, 0) -- (5,1.185);
\draw[black,line width=0.4pt][-] (9.4, 0) -- (9.4,1.185);
\node[align=left] at (2.4, 1.35) {$k=1/R$};
\node[align=left] at (12, 1.35) {$k=1/R$};
\node[align=center] at (5, -0.5) {$\sigma=-1/(\sqrt{7}R)$};
\node[align=center] at (10, -0.5) {$\sigma=1/(\sqrt{7}R)$};
\node[align=left] at (9, 2.2)  {\footnotesize${k=\displaystyle{\frac{8}{7R}} -c\sigma^2}$}; 
\begin{axis}[
    xtick=\empty, 
    ytick=\empty,
  axis lines=middle,   		
  width=16cm,
  height=6cm, 
        xmin=-10, xmax=10,         
        ymin=0, ymax=1,
    xlabel = $\sigma$,
    ylabel = {$k$},
    xticklabel=\empty,
    yticklabel=\empty,
]
\addplot [
line width=1pt,
    domain=-10:10, 
    samples=100, 
    color=black,
   ]
  {0.25+0.1/(1 + (x/2)^4)};    
\end{axis}
\end{tikzpicture}
\caption{The initial curvature profile $k=\mathcal{K}_{1}(\protect\sigma )$
of the generating curve}
\label{fig:RaziPlot2}
\end{figure}
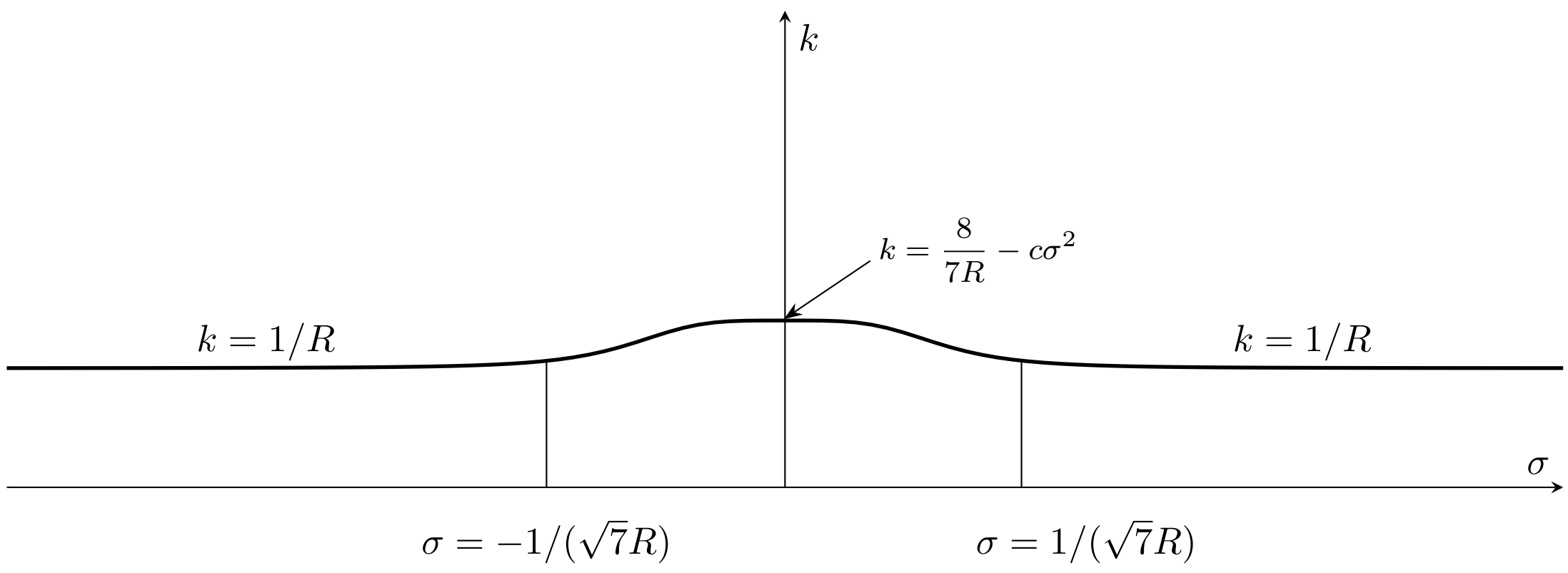

\subsection{The construction of generating curves $\protect\gamma %
_{0}\subset \mathbb{E}^{2}$}

First of all, if $\gamma _{0}\subset \mathbb{E}^{2}$ is the generating curve
of an $O(n)$-invariant , closed $C^{\infty }$-hypersurface $M^{n}\subset 
\mathbb{E}^{n+1}$ with the PRE-property, then $\gamma _{0}$ will be a $%
C^{\infty }$-curve with initial and terminal points perpendicular to the $x$%
-axis and characterized by its \emph{curvature profile}, $k_{1}(\sigma )$,
namely curvature as a function of its arc-length parameter $\sigma \in
\lbrack 0,l(\gamma _{0})]$, up to an $x$-translation. Therefore, for the
purpose of constructing counterexamples as stated above, it is natural and
advantageous to use the setting of Section 4 and try to make use of the
remarkably neat form of the Ricci evolution equation of $k_{1}(\sigma )$ in
order to construct the initial $k_{1}(\sigma _{0})$ so that the \emph{proof
of non-preservation }of positivity in an extremely short time would be very
simple and easy. Note that in the case where the \emph{curvature profile}
with its parameter $\sigma _{0}$ equals $0$ at its symmetry middle point,
the Ricci evolution of $k_{1}(\sigma _{t})$ will also be symmetric with
respect to the $\sigma _{t}\equiv 0$ point and $\lambda _{t}(0)\equiv 0$.
This initial curvature profile will be constructed as follows (see Figure 2).

We start with the $C^{0}$-function $\mathcal{K}_{1}(\sigma )$ defined by%
\begin{equation*}
\mathcal{K}_{1}(\sigma )=\left\{ 
\begin{array}{c}
\frac{8}{7R}-c\sigma ^{2},\text{ \ }|\sigma |\leq \frac{1}{\sqrt{7}}\frac{1}{%
R}, \\ 
\frac{1}{R}\text{, }|\sigma |\geq \frac{1}{\sqrt{7}}\frac{1}{R}, \text{ \ \
\ \ \ \ \ \ \ }R\gg 1%
\end{array}%
\right.
\end{equation*}%
which has corners at $\sigma =\pm \frac{1}{\sqrt{7}}\frac{1}{R}$, and $c
\equiv R \text{ }(mod\frac{1}{R})$. Then one modifies the function with a
rather gentle $C^{\infty }$-smoothing over a considerably wide neighborhood
of $\sigma =\pm \frac{1}{\sqrt{7}}\frac{1}{R}$ to obtain a reflectional
symmetric $C^{\infty }$-curvature profile whose graph around the vicinity of
the origin $\sigma =0$ is as indicated in Figure~\ref{fig:RaziPlot2}. 

Therefore, the Ricci evolution equation (cf. Lemma 3, resp. Lemma 3') at the
middle point , i.e. at $\sigma _{t}=0$, will be simplified as follows,

\begin{equation}
(\mathcal{E}_{3}^{\prime })_{\sigma _{t}=0}:\frac{\partial k_{1}^{t}}{%
\partial t}|_{0}=\frac{\partial ^{2}}{\partial \sigma _{t}^{2}}%
k_{1}^{t}|_{0}-(k_{1}^{t})_{0}^{3}+(n-2)\frac{\partial ^{2}}{\partial \sigma
_{t}^{2}}\left( \frac{\cos ^{2}\lambda }{y}\right) |_{0}  \tag{5.1}
\end{equation}

\noindent and the examples, as illustrated in Figure~\ref{fig:RaziPlot2},
are specifically designed to make the differential calculus of $(\mathcal{E}%
_{2}^{\prime })$ in the vicinity of the short time evolution curve, $\sigma
_{t}=0,\lambda _{t}(0)\equiv 0$, particularly simple to estimate.

One of the basic features of $(\mathcal{E}_{2}^{\prime })$ is that
differentiations of $\frac{\partial }{\partial t}$ can be expressed simply
in terms of that of $\frac{\partial }{\partial \sigma _{t}}$, while 
\begin{equation}
\frac{\partial }{\partial \sigma _{t}}=g(t,\sigma _{0})^{-1}\frac{\partial }{%
\partial \sigma _{0}}\text{ \ \ and \ }\frac{\partial }{\partial t}\ln
g=-2k_{1}k_{2}  \tag{5.2}
\end{equation}%
Thus, in the case of $R\gg 1$ and $T\leq \frac{4}{7}(\frac{1}{R^{2}}+\frac{1%
}{R^{3}})$ one has%
\begin{equation}
k_{1}k_{2}\sim \frac{2}{R^{2}}g(t,\sigma _{0})^{-1}<1+\frac{2t}{R^{2}},\text{
\ \ }0<t<\frac{4}{7}(\frac{1}{R^{2}}+\frac{1}{R^{3}})  \tag{5.3}
\end{equation}%
and therefore $\frac{\partial }{\partial \sigma _{t}}$ is extremely close to 
$\frac{\partial }{\partial \sigma _{0}}$, which makes the task of
estimations of $\frac{\partial ^{l}}{\partial t^{l}}k_{1}^{t}|_{t=0,\sigma
_{0}=0}$ quite straightforward, because in the vicinity \ of $\sigma _{0}$

\begin{equation}
\frac{\partial ^{l}}{\partial \sigma _{0}^{l}}k_{1}^{t}|_{t=0}\equiv 0\text{%
, \ }l\geq 3\text{, \ }\frac{\partial ^{2}}{\partial \sigma _{0}^{2}}%
k_{1}^{t}\equiv -2c\text{ \ }(\mod\frac{1}{R})  \tag{5.4}
\end{equation}

\subsection{Estimation of $(\frac{\partial }{\partial t}%
)^{l}k_{1}^{t}|_{(0,0)}$\ via differential calculus of $(\mathcal{E}%
_{2}^{\prime })$ and $(\mathcal{E}_{3}^{\prime })$ in the vicinity of the
line $\protect\sigma _{t}$\ $\equiv 0,$ $0\leq t$\ $\leq $\ $\frac{a}{R^{2}}$%
\ }

For later calculations or estimates, we introduce the constant $a=0.572$
(e.g. slightly larger than 4/7). Just for the sake of simplifying the
notation of differential calculations in this subsection, we shall introduce
the following system of notation, namely the abbreviations 
\begin{equation}
\left\{ D,\bar{D},d,d_{0}\right\} \text{ \ } \Longleftrightarrow \text{ \ }%
\left\{ \frac{\partial }{\partial \sigma _{t}},\frac{\partial }{\partial
\sigma _{0}},\frac{\partial }{\partial t},\frac{\partial }{\partial t}%
|_{0}\right\}  \tag{5.5} \\
\end{equation}

\begin{equation}
\left\{ k,k_{0},\bar{k},f\right\} \Longleftrightarrow \left\{
k_{1}^{t},k(0),k_{1}^{t}|_{t=0}\ ,\frac{\cos ^{2}\lambda }{y}\right\}  \notag
\end{equation}

Firstly, note that $\bar{k}_{1}(\sigma _{0})$ is an even function, and hence
the Ricci evolution is also reflectionally symmetric, and hence $\cos
\lambda, y$ are even functions while $\sin \lambda ,\tan \lambda $ are odd
functions, and moreover, $D$ applied to even (resp. odd) functions are odd
(resp. even) functions.

Secondly, $D=g(t,\sigma _{0})^{-1}\bar{D}$ and 
\begin{equation}
g(t,\sigma _{0})^{-1}\equiv 1+\frac{2t}{R^{2}}\text{ }(\mod%
\frac{t^{2}}{R^{4}}),\text{ \ }0<t<\frac{4}{7}(\frac{1}{R^{2}}+\frac{1}{R^{3}%
})  \tag{5.6}
\end{equation}%
namely, almost identical up to differences less than $\frac{4t^{2}}{R^{4}}$.
Therefore, by and large, the "non-linearity" of such a differential analysis
of $(\mathcal{E}_{2}^{\prime })$ and $(\mathcal{E}_{3}^{\prime })$ lies in
that of%
\begin{equation}
(n-2)D^{2}(f)\text{ \ }\left( \text{i.e. \ }(n-2)\bar{D}^{2}(f)\right) 
\tag{5.7}
\end{equation}%
Based upon Nash-Moser Theory in the simplest setting of parabolic equations
of $(\mathcal{E}_{0}),(\mathcal{E}_{1})$, and $(\mathcal{E}_{3})$, a concise
simple proof of the \emph{non-preservation } of PRE-property of $%
M^{n}(\gamma _{t})$ is to provide a clean cut \emph{estimate} of $%
k_{1}^{T}(0)$ \emph{less than zero} (cf. Section 1), via Taylor's formula of 
$k_{1}^{t}(0)$, $0\leq t\leq \frac{a}{R^{2}}.$

Therefore, the crucial technicality mainly consists of a highly accurate
estimate for $d_{0}^{l}\bar{k}|_{0}$ of the beginning few $l^{\prime }s$
and, for larger $l^{\prime }s$, simply their order of magnitudes (w.r.t. $%
R^{j}$), together with the specialities of Examples of this section such as
the symmetry and $\bar{D}^{i}\bar{k}\equiv 0$ for $i\geq 3$, in the vicinity
of $\sigma _{0}=0$. For example, the following very simple inductive formula
greatly facilitates the estimation for higher $l^{\prime }s$.

\subsubsection{Simple inductive formulas for estimations of $d_{0}^{l}\bar{k}%
|_{0}$ as well as $d_{0}^{l}\bar{D}^{2}(f)$}

As a direct application of the simplification of $(\mathcal{E}_{3}^{\prime
})|_{\sigma _{t}=0}$, namely%
\begin{equation}
(\mathcal{E}_{3}^{\prime })|_{0}:\frac{\partial k_{1}^{t}}{\partial t}%
|_{\sigma _{t}=0}=\left\{ D^{2}k-k^{3}+(n-2)D^{2}f\right\} |_{\sigma _{t}=0}
\tag{5.8}
\end{equation}%
one has the following simple induction formula:

\begin{equation}
(\circledast )_{l}:\left\{ 
\begin{array}{c}
d_{0}^{(l+1)}k|_{0}=d_{0}^{l}\left\{ D^{2}k-k^{3}+(n-2)D^{2}f\right\} |_{0}%
\text{ \ \ \ \ \ \ \ \ \ \ \ \ \ \ \ \ \ \ \ \ } \\ 
=\left\{ d_{0}^{(l-1)}D^{2}dk-d_{0}^{l}k^{3}+(n-2)d_{0}^{l}D^{2}f\right\}
|_{0}\text{, \ \ }\bar{D}^{4}k|_{0}=0 \\ 
=\left\{
-(D^{2}+d)d_{0}^{(l-1)}k^{3}+(n-2)(D^{2}+d)d_{0}^{(l-1)}(D^{2}f)\right\}
|_{0}%
\end{array}%
\right.  \tag{5.9}
\end{equation}

\noindent \textbf{Remarks 8} \textbf{(i)} Note that $k\sim \frac{2}{R}$ and $%
dk\equiv D^{2}k$ modulo $\frac{1}{R}$ in the vicinity of the segment $\sigma
_{t}=0$, $0\leq t\leq \frac{a}{R^{2}}$. Therefore, intuitively and roughly
speaking, for slightly larger $l$ (say, $l\geq 5$), one should have the
inductive estimate that 
\begin{equation}
d_{o}^{(l+1)}k|_{0}\sim -2d_{0}^{l}(k^{3})|_{0}  \tag{5.10}
\end{equation}%
for the "leading terms" (e.g. order of magnitude).

\textbf{(ii)} The proof of the above inductive estimate amounts to the
estimation of the order of magnitudes of 
\begin{equation*}
(D^{2}+d)d_{0}^{(l-1)}D^{2}f\text{ }
\end{equation*}%
for $(l-1)$ slightly larger than $1$.

\subsubsection{Estimations of the $d_{0}^{l}$ of $k$ (resp. $D^{2}f,D^{4}f)$
for small $l$}

\textbf{(i)} First, recall that $\cos ^{2}\lambda =\frac{1}{2}(1+\cos
2\lambda )$. We calculate the following: \textbf{\ }%
\begin{equation}
D\lambda = -k,\text{ }D\cos ^{2}\lambda =\sin 2\lambda \cdot k  \notag \\
\end{equation}
\begin{equation}
D^{2}\cos ^{2}\lambda = -2\cos 2\lambda \cdot k^{2}+\sin 2\lambda \cdot Dk 
\notag
\end{equation}
\begin{equation}
D^{4}\cos ^{2}\lambda |_{\sigma _{t}=0} =
8k_{0}^{4}-2D^{2}k^{2}|_{0}-4k\cdot D^{2}k|_{0}  \tag{5.11}
\end{equation}
\begin{equation}
\equiv -8kD^{2}k|_{0}=32\frac{c}{R}\text{ }(\mod \frac{1}{R})  \notag
\end{equation}
\begin{equation}
D(y^{-1}) = -y^{-2}\sin \lambda  \notag \\
\end{equation}
\begin{equation}
D^{2}(y^{-1}) = -2y^{-3}\sin ^{2}\lambda +y^{-2}\cos \lambda k\text{, \ }%
D^{3}(y^{-1})|_{0}=0  \notag
\end{equation}
\begin{equation}
D^{4}(y^{-1})|_{0} \equiv y^{-2}\cos \lambda D^{2}k|_{0}=y_{0}^{-2}(-2c)=-%
\frac{2c}{R^{2}}\text{ \ }(mod \frac{1}{R^{2}})  \notag
\end{equation}

Hence,%
\begin{equation}
d_{0}D^{2}f|_{0} = d_{0}\left\{ -\frac{2k^{2}}{y}+\frac{k}{y^{2}}\right\}
|_{0}  \notag
\end{equation}
\begin{equation}
\equiv -4k\frac{d_{0}k}{y}+\frac{d_{0}k}{y^{2}}\equiv 14\frac{c}{R^{2}}\text{
\ } (mod \frac{1}{R^{2}})  \tag{5.12}
\end{equation}
\newline
\begin{equation}
D^{4}f|_{0} \equiv \frac{1}{y}D^{4}\cos ^{2}\lambda
|_{0}+D^{4}(y^{-1})|_{0}\equiv 30\frac{c}{R^{2}}\text{ }(mod \frac{1}{R^{2}})
\notag
\end{equation}
\newline

Thus, it follows from $(\circledast )_{1}$ that%
\begin{equation}
d_{0}^{2}k|_{0} = -(D^{2}+d)(k^{3})|_{0}+(n-2)(D^{4}+dD^{2})f|_{0} 
\tag{5.13}
\end{equation}
\begin{equation}
\equiv \left[ 48+44(n-2)\right] \frac{c}{R^{2}}\text{ \ }(mod \frac{1}{R^{2}}%
)  \notag
\end{equation}

\textbf{(ii)} The estimation of $d_{0}^{3}k|_{0}:$%
\begin{eqnarray}
d_{0}^{2}D^{2}f|_{0} &=&d_{0}^{2}\left\{ -\frac{2k^{2}}{y}+\frac{k}{y^{2}}%
\right\} |_{0}\equiv -\frac{4(d_{0}k)^{2}}{y}|_{0}  \notag \\
&\equiv &-16\frac{c^{2}}{R}\text{ }(mod \frac{1}{R})  \notag
\end{eqnarray}
\begin{equation}
d_{0}D^{4}f|_{0} \equiv d_{0}D^{4}\frac{\cos ^{2}\lambda }{y}|_{0} 
\tag{5.14}
\end{equation}
\begin{eqnarray}
\equiv -8d_{0}\left\{ k\frac{D^{2}k}{y}\right\}|_{0}\equiv -\frac{8}{y}%
d_{0}kD^{2}k|_{0} &\equiv &-32\frac{c^{2}}{R}\text{ }(mod\frac{1}{R})  \notag
\end{eqnarray}

Therefore, from $(\circledast )_{2}$ one has

\begin{equation}
d_{0}^{3}k|_{0} = \left\{
-(D^{2}+d_{0})d_{0}(k^{3})+(n-2)(D^{2}+d_{0})d_{0}D^{2}f\right\} |_{0} 
\tag{5.15}
\end{equation}
\begin{equation}
\equiv -\left[ 96+48(n-2)\right] \frac{c^{2}}{R} =-48n\frac{c^{2}}{R}\text{ }%
, \ c\sim R .  \notag
\end{equation}

\textbf{(iii)} The estimation of $d_{0}^{4}k|_{0}$ and qualitative
estimation of the orders of magnitude of $d_{0}^{l}k|_{0}$ for $l\geq 5$:
First of all, note that it follows from $(\circledast )$:\newline

\qquad \textbf{(a)} The signs of $d_{0}^{l}k|_{0}$ are alternating, i.e. $%
(-1)^{l}$.

\qquad \textbf{(b)} The leading term of $d_{0}^{l}k|_{0}, l\geq 4$, is
actually given by 
\begin{equation}
-(D^{2}+d)d^{(l-2)}k^{3}=-2d^{(l-1)}k^{3}  \tag{5.16}
\end{equation}%
such as 
\begin{equation}
d_{0}^{4}k|_{0} \sim -2d_{0}^{3}k^{3}\sim -12(d_{0}k)^{3}= 96c^{3} 
\tag{5.17}
\end{equation}
\begin{equation}
d_{0}^{5}k|_{0} \sim -2d_{0}^{4}k^{3}\text{, etc.}  \notag
\end{equation}

\noindent Therefore, it is quite straightforward to use the generalized
Leibniz formula to compute such a leading term of $-2d_{0}^{(l-1)}k^{3}$
inductively. However, just for the purpose of providing an upper bound for $%
k(\gamma _{T}(0))$: $T=\frac{4}{7}(\frac{1}{R^{2}}+\frac{1}{R^{3}})$, (cf.
Sect. 5.4.) it suffices to show that%
\begin{equation}
d_{0}^{l}k|_{0}\equiv 0\text{ }(mod\text{ }R^{(l-1)}\text{ \ for }l\geq 4. 
\tag{5.18}
\end{equation}

\noindent which is almost trivial to prove inductively. [ In terms of
Taylor's local expansion of $k(\gamma _{T}(0))$, $T\leq \frac{4}{7}(\frac{1}{%
R^{2}}+\frac{1}{R^{3}})$, it corresponds to $T^{l}\cdot
d_{0}^{l}k|_{0}\equiv 0$ $(mod\frac{1}{R^{l}})$ ]

\subsection{The proof of non-preservation of the PRE-property of those
examples of \S\ 5.2}

Summarizing the estimations of \S\ 5.3 on those coefficients of short term
Ricci evolution of the curvature profile $k_{1}^{t}(\sigma _{t})$ of the
symmetric point $\sigma _{t}=0$, namely

$d_{0}^{l}k|_{0},l=1,2,3$ etc., one has in the case of $c=R$, the following
simplified results:
\begin{eqnarray*}
d_{0}k|_{0} &=&-2R \\
d_{0}^{2}k|_{0} &\equiv &0\text{ }\pmod{\frac{1}{R}}\text{, }(\text{%
i.e. of leading order of }\frac{1}{R}) \\
d_{0}^{3}k|_{0} &\equiv &0\text{ }\pmod {R} \\
d_{0}^{l}k|_{0} &\equiv &0\text{ }\pmod{R^{(l-1)}}\text{ \ for }l\geq 4
\end{eqnarray*}%
while the signs of $d_{0}^{l}k|_{0}$ are alternating, i.e., $(-1)^{l}$. Thus,
the main result of this paper can be simply stated as the following:
\\
\\
\textbf{Main Theorem} We have the estimate
 \begin{center} 
$k(\gamma _{t}(0))\equiv k(\gamma _{0}(0))+d_{0}k|_{0}t\text{    }\pmod{
\frac{1}{R^{5}}}  \text{ \ for }t\leq \frac{4}{7}(\frac{1}{R^{2}}+\frac{1}{%
R^{3}})$,
\end{center}
hence $k(\gamma _{t}(0)$ is no longer positive at $T=\frac{4}{7}(\frac{1%
}{R^{2}}+\frac{1}{R^{3}})$, therefore those examples of \S 5.2 provide simple
counterexamples to the theorem of preservation of \ PRE-property of [4]. 
\\
\\
\textbf{Proof}
It follows directly from the estimate of \S\ 5.3 that all the other terms
of $k(\gamma _{T}(0))$ are at most of the leading orders of $\frac{1}{R^{5}}$%
. In fact, those terms of $\frac{1}{l!}d_{0}^{l}k|_{0}T^{l}$ with $l>5$ are
actually of higher order of smallness, and it is a power series of
alternating signs of decreasing magnitudes.\\

\noindent \textbf{Remark 9}
In retrospect, just for the purpose of proving the main theorem of this
note, it suffices to estimate that the leading orders of $%
d_{0}^{2}k|_{0},d_{0}^{3}k|_{0}$ are, respectively $\frac{1}{R}$ and $R$,
namely, the precise multiples of their leading terms actually play no r%
\^{o}le in the proof.

\section{Concluding remarks}

\qquad \textbf{(1)} Both conceptually and technically, the study of Ricci
evolutions of $O(n)$-hypersurfaces in $\mathbb{E}^{n+1}$ is simply a direct,
straightforward extension of the classical rotational surfaces in $\mathbb{E}%
^{3}$ in terms of the geometric invariant setting of plane curves (cf.
Sections 2, 3, and 4) while the Ricci evolution equation of the fundamental
geometric invariant of plane curves, e.g. the \emph{curvature profile, }$%
k_{1}(\sigma )$, as statedin Lemma 3, concisely and naturally summarizes the
highlight, which turns out to be a rather simple non-linear parabolic
equation with leading term (cf. Eq. (4.14))

\begin{equation}
\frac{\partial }{\partial t}y-\frac{\partial ^{2}}{\partial \sigma _{t}^{2}}%
y,  \tag{6.1}
\end{equation}%
in which the non-linearities are encoded in the time-dependence of $d\sigma
_{t}$ and $\tan \lambda _{t}$, given by the coupled equation:%
\begin{equation}
\frac{\partial g}{\partial t}=-(n-1)k_{1}k_{2}g  \tag{6.2}
\end{equation}%
where $k_{1}k_{2}$ is the Gaussian curvature, i.e. $-\frac{1}{y}\frac{%
\partial ^{2}}{\partial \sigma _{t}^{2}}y$.

\qquad \textbf{(2) }As is often the case of studying the actual solutions as
well as applications of such specific non-linear parabolic PDEs, although
quite simple, it always demands specific geometric-analytic insights and
choice of pertinent applications rather than "general theory". Actually,
this is often where the attractions and challenges lie, while here, as it
turns out, one of the outstanding pertinent applications is the discovery of
the simple family of counterexamples of Section 5.

\qquad \textbf{(3)} As has already been mentioned in Section 5, the specific
geometric-analytic form clearly suggests that the evolution of $k_{1}(\sigma
_{t})$ at the symmetric middle point (i.e. naturally the chosen parameter
origin of $\sigma _{t}=0$) satisfies the simplification of $\lambda
_{t}\equiv 0$ and hence $\tan \lambda _{t}\equiv 0$ . Therefore, the
super-short time analysis of such an almost infinitesimal modification of
the spherical one (i.e. $k_{1}(\sigma )\equiv \frac{1}{R}$, $R\gg 1$)
expectedly becomes very straightforward to compute and to provide clean-cut
estimates with high accuracy. \vspace{\baselineskip}

Geometrically speaking, the remarkable outstanding but intrinsically
spherical (except inside the $\frac{1}{R}$-neighborhood of the geodesic
equator, $S^{n-1}(R)$) all by itself seems to be particularly special, but ,
actually, they have the advantage of \emph{localizability} to that of the $%
\frac{1}{R}$-vicinity of subregions or submanifolds of $S^{n-1}(R)$, which
can be transplanted to all kinds of Riemannian manifolds containing
spherical pieces.

Therefore, for any closed Riemannian manifold $M^{n},n\geq 3$, containing
various locally spherical pieces, each of them can have the above type of
modifications inside of the $\frac{1}{R}$- vicinities of totally geodesic
pieces with all kinds of topological structures. The extremely short Ricci
evolution of such very slightly modified manifolds $\tilde{M}^{n}$ already
contains negative-curvature regions with various versatile topologies. Thus,
the simple construction of $O(n)$-invariant, almost infinitesimal
modifications of $S^{n}(R),R\gg 1$, together with their localizations
already provide a very versatile family of counterexamples towards almost
the proclaimed, celebrated achievements via Ricci evolutions such as [6, 7,
8, 14, 15], namely to the \emph{proofs via Ricci evolutions }but not to the
conclusions or results claimed.

However, the existence of such a family of examples, as it seems to us,
constitutes a clear indication of the following reality of the nature of
Ricci evolutions in general, namely:

\qquad \textbf{(i)} In view of the "basic lemmas" in those works based on
the assumptions of curvature-positivities, e.g. the preservation of such
curvature-positivities are all proved as consequences of Hamilton's "maximal
principle" on Ricci flows; on the one hand, this is the root of faulty
proofs in those papers, and on the other hand, the family of versatile
examples are also counterexamples to the "maximal principle" formulated and
proved via Hamilton's \emph{modified equations }of Ricci flows [4, 5].
Therefore, somewhere and somehow, here must be the origin of
misunderstanding which led to the fallacious developments via Ricci
evolutions through three decades since 1982.

\qquad \textbf{(ii)} A reading of [4, 5] shows that Hamilton's reformulation
and replacement of the Ricci flow equation actually constitute a radical
departure from reality to another system which has radically different basic
analytical features. Thus, the maximal principle he proclaimed is actually
that of his \emph{unrealistic }reformulations which is, in fact, \emph{%
invalid } for Ricci evolution equations.

\qquad \textbf{(iii)}A retrospective reflection on the above examples as
well as reality checks clearly indicates that, on the one hand, the analysis
of the Ricci flow equations is far, far from being adequately understood,
and on the other hand, those simple $O(n)$-hypersurface examples together
with their localizations already shows the proclaimed grand scheme of proof
via Ricci evolutions are simply \emph{misleading }and \emph{unworkable.}

\qquad \textbf{(iv) }Finally, let us take a moment to contemplate on the
perspectives of the global geometry of closed manifolds with curvature
positive property such as the "differentiable sphere problem" as stated in
[8]; and the proof of Poincar\'{e}'s conjecture of dimension 3 or topology
of 3-manifolds, etc. and, moreover, to think about what could be some
reasonable approaches of studying them?

\subsection{Final remarks on the present state of affairs}

\noindent Note that our previous examples only demonstrate that the proofs
via Ricci evolutions of those celebrated works are invalid, and hence these
problems, such as the differentiable sphere problem, 3-dimensional Poincar%
\'{e}'s conjecture etc. remain open and important outstanding problems, and
efforts should be focused on solving them using \emph{other approaches}.

Generally speaking, we think the lesson we learned from the mistakes of
fantasy-chasing of three decades since 1982, based on the misconception of
the Ricci flow equation, due to the lack of basic, realistic understanding
tested by down-to-earth examples, is that one must be \emph{going back to
basics}, grounded on solid, understood knowledge of the topic supported by
realistic testing cases, and patiently approaching via step-wise study of
pertinent intermediate problems.

As an example, let us briefly discuss the "differentiable sphere" as
formulated in [8], but with invalid proof via Ricci flows: The problem,
realistically reformulated, should be what are the best global or localized
pinching bounds for \emph{exotic spheres }of a given dimension? Or more
specifically, what is the best local or global pinching bound for a given or
group of exotic spheres, such as n-dimensional exotic spheres bounding
parallelizable (n+1)-manifolds?

For example, the intent goal of proving upper bound less than $\frac{1}{4}$
is, intuitively, quite far from optimal for exotic spheres, and, on the
other hand, up to now, there are no constructed metrics on exotic spheres of
Milnor or Kervaire type with concrete positive pinching bound.

Historically, Milnor's discovery of 7-dimensional exotic spheres and the
subsequent development of systematic study of \emph{groups of homotopy
spheres }[16] was the grand opening of\emph{\ differential topology}, which
should naturally be followed by efforts to better understand their geometry
and representations as algebraic submanifolds such as the above pinching
problem. Anyhow, it is quite natural to pose as well \ as to pursue the kind
of problems of what are these "best" Riemannian metrics (resp. algebraic
representations) on those exotic spheres of Milnor or Kervaire type with
respect to a naturally chosen property?

The following three kinds of such natural properties have been studied,
namely:

\qquad \textbf{(a)} Using $\dim ISO(\Sigma ,ds^{2})$ as a measurement of its 
\emph{degree of symmetry, }what will be those metrics on a given exotic
sphere $\Sigma $ that will achieve the highest degree of symmetry?

\qquad \textbf{(b) }Defining the \emph{cohomogeneity} of $(\Sigma ,ds^{2})$
to be the dimension of the orbit space $\Sigma /ISO(\Sigma ,ds^{2})$, what
will be those metrics on a given exotic sphere $\Sigma $ that achieve the 
\emph{lowest cohomogeneity?}

\qquad \textbf{(c) }How to represent a given exotic sphere $\Sigma $, say of
Milnor or Kervaire type, as simply as possible as real algebraic sets?%
\newline

It is a pleasant surprise that, not only all of the above three kinds of 
\emph{optimalities} have remarkable clean-cut answers, but they also have
wonderful \emph{concurrence} which reveals the intrinsic harmony between
them, namely the optimalities of (c) are, in fact, isometric representations
of the optimalities of both (a) and (b), while those optimalities of Milnor
spheres (resp. Kervaire spheres) of the same dimension are uniformly equal,
in particular, the cohomogeneity of Kervaire spheres of dimension $(4k+1)$
are all equal to $1$ (cf. [17] and [18]).

On the other hand, the geometry of \emph{cohomogeneity 1 }Riemannian
structures on Kervaire spheres has been studied in [19] via the setting of
equivariant differential geometry of [20] in terms of Killing vector fields
and orbital submersion formula of the Riemannian curvature tensor [21]. In
particular, it is quite straightforward to show that they admit non-negative
sectional curvatures everywhere and in every direction, but necessarily
having zero-sectional curvatures at the singular orbit of codimension $2$,
or in another way.

The above facts can also be seen from equivariant Jacobi-Morse theory on the
global geometry of geodesics. Therefore, the \emph{geometric reality} should
be that those exotic spheres bounding parallelizable manifolds admit metrics
of \emph{non-negative} but \emph{not strictly positive }sectional curvature
everywhere or otherwise. However, it seems to be a technical and rather difficult
proposition to prove (i.e. its intrinsic geometric significance is far
outweighed by its technical difficulty).

\section*{Acknowledgments}

The authors are very grateful to Eldar Straume for his help as well as his pertinent suggestions.
\newline
The conversion of the manuscript prepared by the authors into a \LaTeX\
document featuring embedded TikZ figures was carried out by Eldar Straume
and Razi Naqvi (both at NTNU, Trondheim).

\bigskip

\bigskip

\bigskip

\bigskip

\bigskip


\begin{thebibliography}{99}
\bibitem{} H. E. Rauch. A contribution to differential geometry in the
large. \emph{Annals of Mathematics, }54(1): 38--55, 1951.

\bibitem{} M, Berger. Les vari\'{e}t\'{e}s riemanniennes (1/4)-pinc\'{e}es. 
\emph{Annali della Scuola Normale Superiore di Pisa -- Scienze Fisiche e
Matematiche}, 3e \'{s}erie, 14(2): 161--170, 1960.

\bibitem{} W. Klingenberg. \"{U}ber Riemannsche Mannigfaltigkeiten mit
positiver Kr\"{u}mmung\emph{. Commentarii mathematici Helvetici, }35:
47--54, 1961.

\bibitem{} R. S. Hamilton. Three-manifolds with positive Ricci curvature. 
\emph{Journal of Differential Geometry. }17(2): 255--306, January 1982.

\bibitem{} R. S. Hamilton. Four-manifolds with positive curvature operator. 
\emph{Journal of Differential Geometry} 24(2): 153--179, 1986.

\bibitem{} Ben Andrews and Huy Nguyen. Four-manifolds with $1/4$-pinched
flag curvatures. \emph{Asian Journal of Mathematics, }13(2): 251--270, 2009.

\bibitem{} C. B\"{o}hm and B. Wilking. Manifolds with positive curvature
operators are space forms. \emph{Annals of Mathematics, }167(3): 1079--1097,
2008.

\bibitem{} S. Brendle and R. Schoen. Manifolds with $1/4$-pinched curvature
are space forms. \emph{J. Amer. Math. Soc., }22(1): 287--307, 2009.

\bibitem{} Haiwen Chen. Pointwise $1/4$-pinched 4-manifolds. \emph{Annals of
Global Analysis and Geometry,} 9(2): 161--176, June 1991.

\bibitem{} G. Huisken. Ricci deformation of the metric on a Riemannian
manifold. \emph{Journal of Differential Geometry, }21(1): 47--62, 1985.

\bibitem{} M. J. Micallef and J. Douglas Moore. Minimal two-spheres and the
topology o f manifolds with positive curvature on totally isotropic
two-planes. \emph{Annals of Mathematics, }127(1): 199--227, 1998.

\bibitem{} Huy Nguyen. \emph{Positive isotropic curvature and Ricci flow. }%
PhD thesis, 2008.

\bibitem{} S. Nishikawa. Deformation of Riemannian metrics and manifolds
with bounded curvature ratios. In \emph{Geometric measure theory and the
calculus of variations (Arcata, Calif., 1984), }volume 44 of \emph{Proc.
Sympos. Pure Math., }pages 343--352. Amer. Math. Soc., Providence, RI, 1986.

\bibitem{} G. Perelman. \emph{The entropy formula for the Ricci flow and its
geometric applications.} arXiv: math/0211159, 11 Nov 2002.

\bibitem{} G. Perelman. \emph{\ Ricci flow with surgery on three-manifolds. }%
arXiv: math/0303109, 10 Mar 2003.

\bibitem{} M. A. Kervaire and J. W. Milnor. Groups of homotopy spheres: I. 
\emph{Annals of Mathematics, }77(3): 504--537, 1963.

\bibitem{} E. V. Brieskorn. Examples of singular normal complex spaces which
are topological manifolds. \emph{Proceedings of the National Academy of
Sciences, }55(6): 1395--1397, 1996.

\bibitem{} J. W. Milnor. \emph{Singular Points of Complex Hypersurfaces.
(AM-61). }Princeton University Press, 1968.

\bibitem{} A. Back and Wu-Yi Hsiang. Equivariant geometry and Kervaire
spheres. \emph{Transactions of the American Mathematical Society, }304(1):
207--227, 1987.

\bibitem{} A. Back, M. P. do Carmo, and Wu-Yi Hsiang. On some fundamental
equations of equivariant \ Riemannian geometry. \emph{Tamkang Journal of
Mathematics, }40(4): 343--376, Dec. 2009.

\bibitem{} B. O'Neill. The fundamental equations of a submersion. \emph{%
Michigan Mathematical Journal, }13(4): 459--469, 1966.
\end{thebibliography}
\end{document}